\documentclass{article}%
\usepackage{amsfonts}
\usepackage{amsmath}
\usepackage{amssymb}
\usepackage{color}

\newtheorem{theorem}{Theorem}
\newtheorem{acknowledgement}{Acknowledgement}

\newtheorem{definition}{Definition}
\newtheorem{lemma}{Lemma}

\newtheorem{remark}{Remark}

\begin{document}
\begin{center}
{\Large {\bf On the Spherical Hausdorff Measure  in Step 2 Corank 2 sub-Riemannian Geometry}}

\vskip .5cm
{Ugo Boscain}

{\footnotesize CNRS, CMAP, Ecole Polytechnique, Route de Saclay, 91128 Palaiseau Cedex, and
INRIA Saclay team GECO
{\tt boscain@cmap.polytechnique.fr }}

\vskip .5cm
{Jean-Paul Gauthier}

{\footnotesize LSIS, UMR\ CNRS 6168, Avenue Escadrille Normandie Niemen,13397 MARSEILLE Cedex
20, and INRIA Saclay Team GECO
{\tt gauthier@univ-tln.fr}}

\vskip 1cm
{Dedicated to A.\ Agrachev for his 60th birthday}

\end{center}

\vskip 1cm
\begin{abstract}
In this paper, we \ consider generic corank 2 {sub-Riemannian} structures, and we
show that the Spherical Hausdorf measure is always a $\mathcal{C}^{1}$-smooth
volume, which is in fact generically $\mathcal{C}^{2}$-smooth out of a
stratified subset of codimension 7. In particular, for rank 4, it is
generically $\mathcal{C}^{2}.$
This is  the continuation of a  previous work by the auhors.

\end{abstract}

\vskip .5cm
\noindent
{\bf subjclass: 53C17, 49J15, 58C35}

\section{ \ Introduction}

In this paper we consider sub-Riemannian structures $s=(\Delta,g)$ over an  oriented
$n$-dimensional manifold $M$.\ The distribution $\Delta$ has rank $p$ and
corank $k=n-p,\ $\ and $g$ is a Riemannian metric over $\Delta.\ $In most of
the paper, $k=2.$ Moreover, the distribution is assumed to be 2-step bracket generating.

The set $\mathcal{S}$ of such (corank 2, 2-step bracket generating)
 sub-Riemannian structures over $M$ is  endowed with the $\mathcal{C}^{\infty}$
Whitney topology.

As it will be recalled in the next section, there is a natural smooth measure
associated with the structure $s$, called the Popp measure (see \cite{mont}).
It has been shown in \cite{AB} that the Radon-Nykodim derivative 
 $f_{\cal S P}(\xi)$
of the spherical Haussdorf measure with respect to the Popp measure at a point $\xi$ is just  (universally)
proportional to the inverse of the Popp-volume of the unit ball of the nilpotent approximation of
$s$ at $\xi$. Moreover, in the same paper, when $k=1,$ it is shown
that  $f_{\cal S P}(\xi)$ is a $\mathcal{C}^{3}$ function  ($\mathcal{C}^{4}$ along curves), which is not $\mathcal{C}%
^{5}$ in general.

 The nonsmoothness appears only in sub-Riemannian structures for which the nilpotent approximation depends on the point and can occur at points (called resonance points) where certain invariants of the structures coincide.
The  (high) degree 3 of differentiability is due to the fact that, in the corank 1 case,
the conjugate locus of the nilpotent approximation coincides with the cut
locus.\ This  coincidence is no more true for higher corank. In particular, this is shown in \cite{BBG}, in the corank 2 case, and an
explicit characterization of the cut-locus is given. In the same paper,  as a simple byproduct, it has
been shown that  $f_{\cal S P}$ is generically $\mathcal{C}^{1}$ for ${p=4},k=2.$

\bigskip

Starting from the explicit characterization of the cut-locus obtained in \cite{BBG}, in the current paper we go one step further and obtain the following result.

\begin{theorem}
\label{th} (step 2, corank $2$) We have the following: 

1. the Radon-Nykodim derivative  $f_{\cal S P}$ is
always $\mathcal{C}^{1}$;

2.\ The Radon-Nykodim derivative  $f_{\cal S P}$ is generically\footnote{In this
theorem, genericity means that the property is satisfied for a subset of
 sub-Riemannian metrics that CONTAINS an open-dense set. Indeed, in the
transversality arguments, we can always avoid the closure of certain
Whitney-stratified bad-sets, in place of avoiding the bad-sets themselves.}
$\mathcal{C}^{2},$ out of a stratified set of codimension 7.
In the particular
case {$p=4,$} there is an open-dense subset of $\mathcal{S}$ for which
 $f_{\cal S P}$ is $\mathcal{C}^{2}$-smooth.
\end{theorem}

\begin{remark}
 In the case of a non-orientable manifold the Popp measure cannot be defined as a volume form, but just as a density. However, Theorem \ref{th} still holds true since it is essentially local.   
\end{remark}

Roughly speaking, $f_{\cal S P}$ depends on the maximum eigenvalue of a certain skew symmetric matrix (depending on the point) defining the nilpotent approximation of the structure at the point. (This eigenvalue is an invariant of the structure.) Hence, the study of the differentiability properties of  $f_{\cal S P}$ requires a fine analysis of the regularity of the maximum eigenvalue of a family of skew symmetric matrices smoothly depending on parameters. When the maximum eigenvalue is simple at a point, then in a neighborhood of that point it is ${\cal C}^\infty$. A drop of regularity appears at points where the maximum eigenvalue is multiple. The ${\cal C}^1$ regularity can be obtained as a consequence of the fact that when  the maximum eigenvalue is multiple, the cut locus coincides with the conjugate locus. \footnote{Notice that in the corank 1 case, the cut locus coincides with the conjugate locus at every point, see \cite{AB}.}
This fact does not permit to get  the ${\cal C}^2$  result, which requires a deeper analysis.  To treat double eigenvalues we need an adaptation of a  deep result of Arnold \cite{Arnold}, to the case of versal deformations of real skew-symmetric matrices. The case of triple eigenvalues is apparently extremely difficult and we do not treat it in this paper. 
However,  the set of skew symmetric matrices with a triple eigenvalue is an algebraic subset of codimension 8 in skew symmetric matrices (we provide a proof of this technical fact in appendix). In the particular case of rank 4 this set is generically empty.

The paper is organized as follows: in  Section \ref{next}, we recall
the definition of the Popp measure and that of the nilpotent approximation of
$s=(\Delta,g).$  We avoid to recall all standard definitions of sub-Riemannian geometry since these are already given in  \cite{BBG}. 
Then, we recall the main result of \cite{BBG} which gives the cut time for geodesics issued from the origin.
This is our key point. In Section \ref{proof}, we give the proof of Theorem \ref{th}. 
This proof uses a certain number of technical tools that are collected in appendix.
In Appendix \ref{quat} we recall certain basic facts about
quaternions, which here represents a very convenient tool. In \ref{versal} we study versal deformations of real skew-symmetric matrices. In \ref{a-codimension} we discuss the codimension of  the set of skew symmetric matrices having either a double or a  triple eigenvalue. In \ref{genericityR} we prove a result that  (generically) allow us to make a crucial change of coordinates. In \ref{formula}, we recall how to get a useful formula for the volume of the nilpotent ball in the corank 2 case.

\section{Prerequisites\label{next}}

\subsection{\bigskip Nilpotent approximation\label{nilap}}

We define the nilpotent approximation in the two-step baracket generating case
only.\ The tensor mapping:
\begin{equation}
\lbrack.,.]:\Delta_{\xi}\times\Delta_{\xi}\rightarrow T_{\xi}M/\Delta_{\xi},
\label{brac}%
\end{equation}
is skew symmetric.\ Then, for any $Z^{\ast}\in(T_{\xi}M/\Delta_{\xi})^{\ast}$
we have:%
\[
Z^{\ast}([X,Y]+\Delta_{\xi})=<A_{Z^{\ast}}(X),Y>_{g}%
\]
for some $g$-skew-symmetric endomorphism $A_{Z^{\ast}}$ of $\Delta_{\xi}.$ The
mapping $Z^{\ast}\rightarrow A_{Z^{\ast}}$ is linear, and its image is denoted
by $\mathcal{L}_{\xi}.$

The space $L_{\xi}=\Delta_{\xi}\oplus T_{\xi}M/\Delta_{\xi}$ is endowed with
the structure of a 2-step nilpotent Lie algebra with the bracket:
\[
\lbrack(V_{1},W_{1}),(V_{2},W_{2})]=(0,[V_{1},V_{2}]+\Delta_{q}).
\]
The associated simply connected nilpotent Lie group is denoted by $G_{\xi},$
and the exponential mapping $E_{xp}:L_{\xi}\rightarrow G_{\xi}$ is one-to-one
and onto. By translation, the metric $g_{\xi}$ over $\Delta_{\xi}$ allows to
define a left-invariant  sub-Riemannian structure over $G_{\xi},$ called the
nilpotent approximation of $(\Delta,g)$ at $\xi.$

Any $k$-dimensional vector subspace $\mathcal{V}_{\xi}$ of $T_{\xi}M$,
transversal to $\Delta_{\xi},$ allows to identify $L_{\xi}$ and $G_{\xi}$ to
$T_{\xi}M\simeq\Delta_{\xi}\oplus T_{\xi}M/\Delta_{\xi}.$ If we fix $\xi
_{0}\in M,$we can chose linear coordinates $x$ in $\Delta_{\xi_{0}}$ such that
the metric $g_{\xi_{0}}$ is the standard Euclidean metric, and for any linear
coordinate system $y$ in $\mathcal{V}_{\xi_0},$ there are skew-symmetric
matrices $L_{1},...,L_{k}$ $\in so(p,\mathbb{R})$ such that te mapping
\ref{brac} writes:%
\[
\lbrack X,Y]+\Delta_{\xi}=\left(
\begin{array}
[c]{c}%
X^{\prime}L_{1}Y\\
.\\
.\\
X^{\prime}L_{k}Y
\end{array}
\right)  ,
\]
where $X^{\prime}$ denotes the transpose of the vector $X.$

This construction works for any $\Delta,$ but $\Delta$ is  2-step bracket
generating iff the endomorphisms of $\Delta_{\xi},$  $L_i$, $i=1,...,k$
(respectively the matices $L_{i}$ if coordinates $y$ in $\mathcal{V}_{\xi}$
are chosen) are linearly independant.

\bigskip

\subsection{Popp Measure\label{popp}}

In the 2-step bracket generating case, the linear coordinates $y$ in $T_{\xi
}M/\Delta_{\xi}$ can be chosen in such a way that the endomorphisms  $L_i$,
$i=1,...,k$ are orthonormal with respect to the Hilbert-Schmidt scalar product
$<L_{i},L_{j}>=\frac{1}{p}\mbox{Trace}_{g}(L_{i}^{\prime}L_{j}).$ This choice defines
a canonical euclidean structure over $T_{\xi}M/\Delta_{\xi}$ and a
corresponding volume in this space.\ Then using the Euclidean structure over
$\Delta_{\xi},$ we get a canonical eucildean structure over $\Delta_{\xi
}\oplus T_{\xi}M/\Delta_{\xi}.\ $The choice of the subspace $\mathcal{V}_{\xi
}$ induces an euclidean structure on $T_{\xi}M$ that depends on the choice of
$\mathcal{V}_{\xi},$ but the associated volume over $T_{\xi}M$ is independant
of this choice.

\begin{definition}
This volume form on $M$ is called the Popp measure.
\end{definition}
 By construction, the Popp measure is a smooth volume form.

\bigskip

Let us recall a main result from \cite{AB}.

\begin{theorem}
\label{thhaus}  (equiregular, any step, any corank) The value 
 $f_{\cal S P}(\xi)$ at $\xi\in M$ of the Radon-Nykodim
derivative of the spherical Hausdorf measure with respect to the Popp measure
is equal to  $2^Q/\hat\mu(\hat B_\xi)$, where $Q$ is the Hausdorff dimension of the sub-Riemannian structure as metric space and $\hat\mu(\hat B_\xi)$ is the Popp volume of the unit ball of the nilpotent approximation at
$\xi.$
\end{theorem}

\subsection{Geodesics and Cut-locus\label{cutl}}

We restrict to the corank 2 case.\ Here, we consider geodesics of the
nilpotent approximation of $s=(\Delta,g)$ in $T_{\xi_{0}}M\simeq\mathbb{R}%
^{n},\ $issued from the origin.\ A transversal subspace $\mathcal{V}_{\xi_{0}%
}$ is chosen, together with the linear Hilbert-Schmidt-orthonormal coordinates
$y$ in $\mathcal{V}_{\xi_{0}}$, and euclidean coordinates $x$ in $\Delta
_{\xi_{0}}.$ The geodesics are projections on $\mathbb{R}^{n}$ of trajectories
of the  smooth  Hamiltonian $H$ on $T^{\ast}\mathbb{R}^{n}:$%

\begin{equation}
H(p^x,p^y,x,y)=\sup_{u\in\mathbb{R}^{p}}(-\|u\|^{2}+%
{\displaystyle\sum\limits_{i=1}^{p}}
p^x_{i}u_{i}+p^y_{1}x^{\prime}L_{1}u+p^y_{2}x^{\prime}L_{2}u). \label{ham}%
\end{equation}
where $p^x,p^y$ are the coordinates dual to $x,y.$ Geodesics are arclength-parametrized as soon as the initial covector 
$(p^x(0),p^y(0))$
verifies $H(p^x(0),p^y(0),x(0),y(0))=1/2$. For geodesics issued from the origin, this condition reads $\|u(0)\|= \|p^x(0)\|=1$, where the norm $\|\cdot\|$ is the one induced by duality on $\Delta_{\xi_0}^\ast$.

Note that $p^y_{1},p^y_{2}$ are constant along geodesics, since the  Hamiltonian
(\ref{ham}) does not depend on the $y$-coordinates. 

The following result is shown in \cite{BBG}, and is crucial for the proof of
our result.

\begin{theorem}
\label{thcut}The cut time $t_{cut}$ \ of the arclength-parametrized geodesic  issued from the origin and
corresponding to  the initial covectors $(p^x(0),p^y)$ is given by:%

\[
 t_{cut}=
\frac{2\pi}{\max(\sigma(p^y_{1}L_{1}+p^y_{2}L_{2})},
\]
where $\max(\sigma(A))$ denotes the maximum modulus of the eigenvalues of the
skew symmetric matrix $A.$  In general, the conjugate time is not equal to the
cut time. \end{theorem}

\begin{remark}
 In fact the cut time is also conjugate if and only if the  matrix $p^y_{1}L_{1}+p^y_{2}L_{2}$ has a double maximum eigenvalue or 
$[L_1,L_2]=0$.
\end{remark}

It turns out that the singularities of the Hausdorf measure appear due to collision between the two largest moduli of eigenvalues of the matrix $p^y_{1}L_{1}+p^y_{2}L_{2}.$ The set of
skew-symmetric matrices that have a double eigenvalue is a codimension 3
algebraic subset of $so(p,\mathbb{R})$ (see Appendix \ref{a-codimension}). Then, from the tranversality theorems
(\cite{abra}), for generic (open, dense)  sub-Riemannian structures, the set
$\overline{\Sigma_{2}}$ of points of $M$ such that $p^y_{1}L_{1}+p^y_{2}L_{2}$ has
a double (at least) eigenvalue for some $p^y_{1},p^y_{2}$ has codimension 2 in
$M.$ The problems of smoothness of the Hausdorf measure will occur on
$\overline{\Sigma_{2}}$ only.

 Along the paper we set,
for the geodesic under consideration:%
\[
p^y_{1}=r\cos(\theta),\text{ }p^y_{2}=r\sin(\theta),\text{ }A_{\xi}(\theta
,r)=\frac{2\pi}{\max(\sigma(p^y_{1}L_{1}+p^y_{2}L_{2} ))},\text{ where }\xi=(x,y)\in
M.
\]

It is known (\cite{bron, kur, milch}) that $A_{\xi}(\theta,r)$ is a Lipschitz
function of all parameters $\xi,\theta,r.$ We write also $A_{\xi}%
(\theta)=A_{\xi}(\theta,1).$

\section{Proof of Theorem \ref{th}}\label{proof}

For a fixed point $\xi_{0}=(x_{0},y_{0})\in M,$ let us consider the
exponential mapping $\mathcal{E}$ associated with the nilpotent approximation
at $\xi_{0}$, where $x,y$ are coordinates as in Section \ref{popp}:%
\[
\mathcal{E}_{t}\mathcal{(}p^x_{0},p^y_{0})=\pi(e^{t\vec{H}}(p^x_{0}
,p^y_{0},{\xi_{0}})),
\]
where $\pi:T^{\ast}M\rightarrow M$ is the canonical projection, and $\vec{H}$
is the  Hamiltonian vector field associated with the  Hamiltonian
(\ref{ham}).  {Here 
$(p^x_0,p^y_0)$ are initial covectors satisfying $H(p^x_0,p^y_0,\xi_0)=1/2$.}

As above we have $p^y(t)=p^y_{0}=(p^y_1,p^y_2)=(r\cos
(\theta),r\sin(\theta))$. Also, by  homogeneity, $\mathcal{E}%
_{t}\mathcal{(} p^x_0,p^y_{0})=\mathcal{E}_{1}\mathcal{(}t$  $p^x_{0},t$ $p^y_{0}).\ $

In our paper \cite{BBG}, the following formula is given for the volume $V_{\xi}$ at a point $\xi\in M$\ of the unit ball
of the nilpotent approximation. For the benefit of the reader, this formula is established here in Appendix \ref{formula}.

\begin{equation}
V_{\xi}=\int_{0}^{2\pi}\int_{0}^{A_{\xi}(\theta)}\int_{B}J( p_0^x,\theta,r,\xi)d p_0^x\text{ }dr\text{ }d\theta\label{vol}%
\end{equation}
where $B$ is the unit ball in the euclidean $p$-dimensional
$ p_0^x$-space, and $J(p_0^x,\theta,r,\xi)$ is the jacobian determinant
of $\mathcal{E}_{1}\mathcal{(}{ p_0^x},r\cos(\theta),r\sin(\theta)).$

We set $f_{\xi}(\theta,r)=\int_{B}J(p_0^x,\theta,r,\xi)d p_0^x,$ and
$W_{\xi}(\theta)=$ $\int_{0}^{A_{\xi}(\theta)}f_{\xi}(\theta,r)dr.$ If we show
that $W_{\xi}(\theta)$ is $\mathcal{C}^{1}$ or $\mathcal{C}^{2}$ w.r.t
$(\theta,\xi),$ it will imply that $V_{\xi}$ is $\mathcal{C}^{1}$ or
$\mathcal{C}^{2}$ w.r.t $\xi.$

 In a neighborhood of a fixed  $(\theta_{0},\xi
_{0})\in  S^1\times M$ we have,
\begin{align}
W_{\xi}(\theta)  &  =\int_{0}^{A_{\xi}(\theta)}f_{\xi}(\theta,r)dr\text{
}\label{wol}\\
&  =\int_{0}^{A_{\xi_{0}}(\theta_{0})}f_{\xi}(\theta,r)dr\text{ }+\int%
_{A_{\xi_{0}}(\theta_{0})}^{A_{\xi}(\theta)}f_{\xi}(\theta,r)dr\nonumber\\
&   =:(I)+(II).\text{ }\nonumber
\end{align}
The term (I) is smooth.\  We are then left to study the smoothness
of $II(\xi,\theta)$.

\subsection{Proof of the fact that $W_{\xi}(\theta)$ is always $\mathcal{C}%
^{1}$\label{prc1}}

Setting $z=(\theta,\xi),$  $z_0=(\theta_0,\xi_0),$ and $f(z,r)=f_{\xi}(\theta,r),$ $A(z)=A_{\xi}%
(\theta)$, the tangent mapping to $II(\xi,\theta),$ at $(\theta_{0},\xi_{0})$ is
\begin{equation}
D\text{ }II(z_{0})(h)=\sum_{i=1}^{n+1}f(z_{0},A(z_{0}))\frac{\partial
A}{\partial z_{i}}(z_{0})h_{i}. \label{firstder}%
\end{equation}

This last expression makes sense, and is continuous w.r.t $z_{0}$ for the
following reasons: first as we said, $A(z)$ is Lipschitz-continuous, then the
derivatives are bounded.  Moreover at points $z_{0}$ such that $A$ is
not differentiable, $f(z_{0},A(z_{0}))$ vanishes.\ This last point follows from
the fact that when the eigenvalue of $A(z_{0})$ having maximum modulus is multiple then
 the conjugate time is
equal to the cut time, which makes the jacobian determinant $J( p_0^x,\theta_{0},A(\theta_{0},\xi_{0}),\xi_{0})$ vanish for all $ p_0^x$. This comes
from the section II.3\ 1 in the paper \cite{AB}.

\bigskip

\begin{remark}
\label{rem}In fact, it follows from the same paper that, if $A(z_{0})$  corresponds to a
multiple eigenvalue, then the rank of 
$J_{\xi_0}( p_0^x,\theta_{0},A(\theta_{0},
\xi_{0}))= J( p_0^x,\theta_{0},A(\theta_{0},
\xi_{0}),\xi_{0})$ drops by 2 at least, independently of  $p_0^x.$ This point will
be very important in the next section.
\end{remark}

This ends the proof.

\subsection{Proof of the $\mathcal{C}^{2}$ result\label{prc2}}

It follows from the transversality theorems (\cite{abra, gor}) and from Lemma
\ref{l1} and Lemma \ref{sub}  in the Appendix, that there exists an open dense  subset of  sub-Riemannian metrics,
still denoted by $\mathcal{S},$ such that all elements $s$ of $\mathcal{S}$
meet: the set $U_{s}\subset  S^1\times M$ of  $(\theta,\xi)$ such that
$A_{\xi}(\theta)$ corresponds to a triple (at least) eigenvalue is a locally
finite union of manifolds, regularly embedded, of codimension 8 in $S^1\times M,$ and the set $\tilde{U}_{s}$ $\subset  S^1\times M$ of $(\theta,\xi)$
such that $A(\theta,\xi)$ corresponds to a double (and not triple) eigenvalue
is a locally finite union of manifolds, of codimension 3.\bigskip

We want to show the following property (P), for a  (smaller) generic (residual in the
Whitney topology) set $\mathcal{S}_{0}$ of  sub-Riemannian metrics over $M$:

\bigskip
(P) the partial derivatives $D_{i}(z)=$ $f(z,A(z))\frac{\partial A}{\partial
z_{i}}(z)$ from (\ref{firstder}) are $\mathcal{C}^{1}$ in a neighborhood of all
points $z_{0}$ such that $A(z_{0})$ corresponds to a double (and not triple) eigenvalue.

\bigskip

To do this, we fix $s_{0}$ and $z_{0}\in\tilde{U}_{s_{0}}$ and we consider a
(mini)versal deformation of $L(\xi_{0}, \theta_0)=L_{1}(\xi_{0})\cos(\theta
_{0})+L_{2}(\xi_{0})\sin(\theta_{0})=L(z_{0}),$ as introduced in  Appendix
\ref{versal}.\ It follows that:%

\[
L(\xi,\theta)=L(z)=g(z)^{-1}
{\cal T}(z)
g(z)
\]
 where  $g(z)$ belongs to the
orthogonal group 
and\footnote{${\cal T }(z)$ plays the role of $T(\mu(z))$ in Appendix \ref{versal}.
${\cal T}(z)$ is the block-diagonal matrix  ${\cal T }(z)=\mbox{Bd}(\lambda(z)\hat{q}+q(z),\Delta(z)).$ Here, following the notation introduced in the appendix, $q$ is a pure quaternion, $\hat q$ is a pure skew-quaternion, $\Delta(z)$ is a $2\times2$ block diagonal skew-symmetric matrix and $\lambda(z)$ is a nonzero real number.
} the functions $g(\cdot),\lambda(\cdot),q(\cdot),\Delta(\cdot)$ are smooth with respect to  $z$.\ 

\bigskip
\noindent 
{The following crucial Lemma is proved in Appendix \ref{genericityR}

\begin{lemma}\label{l-genR}
The property\\[2mm]
(R): the map $ S^1\times M\ni z \mapsto q(z)\in \mathbb{R}^{3},$
has rank 3 at every $z\in \tilde{U}_{s}$,\\[2mm]
is residual in $\mathcal{S}$.
\end{lemma}

Let us call ${\cal S}_0$ the subset of ${\cal S}$ for which (R) holds. If $s_0$ is fixed in ${\cal S}_0$  and $z_{0}\in\tilde{U}_{s_{0}}$  then, locally around $z_0$,
we can find a system of coordinates in $ S^1\times M$ in such a way that the three first
coordinates, $z_{1},z_{2},z_{3}$ become the three components of $q(z)$.
 Note that these 3 coordinates vanish at
$z_{0}.$}

Locally, the codimension 3 manifold $\tilde{U}_{s_{0}}$ is determined by the
equations $z_{1}=z_{2}=z_{3}=0.$

\bigskip

As we said in Remark \ref{rem}, the rank of $J( p_0^x,z,A(z))$ drops
by 2 at least, independantly of $p_0^x,$ at  each point $z\in\tilde{U}_{s_{0}}.$
Formula (\ref{eigquat}) in the appendix tell us that
$A(z)=\frac{2\pi}{\lambda(z)+\sqrt{z_{1}^{2}+z_{2}^{2}+z_{3}^{2}}}$ where
$\lambda(z)$ is smooth and nonzero. We set $\hat{z}_{4}=(z_{4},...,z_{n+1})$
and $\hat{z}_{1}=(z_{1},z_{2},z_{3}).$ 

The Jacobian determinant $J(p_0^x,z,r)$ can be written as
$$
V_1(p_0^x,z,r)\wedge\ldots\wedge V_{n+1}(p_0^x,z,r),
$$
for certain smooth $n+1$-dimensional vectors $V_1(p_0^x,z,r)\ldots V_{n+1}(p_0^x,z,r)$.

For all $p_0^x$, at points $(z,r)$ such that $\hat z_1=0$, $r=A(z)=\frac{2\pi}{\lambda(z)}$, the vectors $V_1\ldots V_{n+1}$ have rank $n-1$ at most. Then

$$
0=\frac{\partial J}{\partial z_i}=\frac{\partial V_1}{\partial z_i}\wedge V_2\wedge\ldots \wedge V_{n+1}+
V_1\wedge\frac{\partial V_2}{\partial z_1}\wedge\ldots \wedge V_{n+1}+ \ldots
$$
and 
$$
0=\frac{\partial J}{\partial r}=\frac{\partial V_1}{\partial r}\wedge V_2\wedge\ldots \wedge V_{n+1}+
V_1\wedge\frac{\partial V_2}{\partial r}\wedge\ldots \wedge V_{n+1}+ \ldots
$$
It follows that $J$, $\frac{\partial J}{\partial z_i}$, $\frac{\partial J}{\partial r}$ vanish at all $(z,r)$ with $\hat z_1=0$, $r=A(z)=\frac{2\pi}{\lambda(z)}$.

Therefore $f(z,r)=\int_{B^1} J_\xi (p_0^x,\theta,r)\,dp_0^x$ is a quadratic expression in the variable $\hat z_1$, $r-\frac{2\pi}{\lambda(z)}$ depending smoothly on $z$, $r$:


\begin{equation}
f(z,r)=\tilde{Q}_{z,r}(\hat{z}_{1},r-\frac{2\pi}{\lambda(z)}). \label{expf}%
\end{equation}

Now we study the continuity of  the second partial derivatives of $W_{\xi}(\theta)=\int%
_{0}^{A(\theta_{0},\xi_{0})}f_{\xi}(\theta,r)dr$ $+\int_{A(\theta_{0},\xi
_{0})}^{A(\theta,\xi)}f_{\xi}(\theta,r)dr,$ or with the new notations,
$W(z)=\int_{0}^{A(z_{0})}f(z,r)dr$ $+\int_{A(z_{0})}^{A(z)}f(z,r)dr.$

The first partial derivatives, at any point $z_{0}$ were:%
\begin{align*}
\frac{\partial W}{\partial z_{i}} (z_{0}) &  =\int_{0}^{A(z_{0})}%
\frac{\partial}{\partial z_{i}}f( z_0,r)dr+f(z_{0},A(z_{0}))\frac{\partial
A}{\partial z_{i}}(z_{0}),\\
&   =:III(z_{0})+IV(z_{0})
\end{align*}

To show that $\frac{\partial III(z)}{\partial zj}$ exists and is continuous,
we proceed exactly as in Section \ref{prc1}, using the fact that$\frac
{\partial}{\partial z_{j}}f(z,r)$ also vanishes at $(\hat{z}_{1}=0$,
$r=\frac{2\pi}{\lambda(z)}).$

The more difficult point is to show that $\frac{\partial IV(z)}{\partial zj}$
exists and is continuous.%
\[
\frac{\partial IV(z)}{\partial  z_j}=\frac{\partial}{\partial z_{j}%
}\left(f(z,A(z))\frac{\partial A(z)}{\partial z_{i}}\right).
\]

We get:%

\begin{align*}
\frac{\partial IV}{\partial  z_j}(z)  &  =\frac{\partial f}{\partial z_{j}%
}(z,A(z))\frac{\partial A(z)}{\partial z_{i}})+\frac{\partial f}{\partial
r}(z,A(z))\frac{\partial A(z)}{\partial z_{i}}\frac{\partial A(z)}{\partial
z_{j}}+f(z,A(z))\frac{\partial^{2}A(z)}{\partial z_{i}\partial z_{j}}.\\
&   =:V(z)+VI(z)+VII(z).
\end{align*}

The cases of $V(z),VI(z)$ are obvious, since again $\frac{\partial
A(z)}{\partial z_{i}}$ is bounded, and the functions $\frac{\partial
f}{\partial z_{j}}(z,A(z)),\frac{\partial f}{\partial r}(z,A(z))$ are
continuous and go to zero when $\hat{z}_{1}$ tends to zero. The only
difficulty is the case of $VII(z).$

Remind that $A(z)=\frac{2\pi}{\lambda(z)+||\hat{z}_{1}||.}$ where $\lambda(z)$
is nonzero, smooth.\ Then the only problem may occur for $i=1,2,3.$

Let us consider only the 2 cases:  (1) $i=1, j=4,$ (2) $i=1, j=2,$ the other 
being similar.

Case (1): $\frac{\partial A(z)}{\partial z_{1}}=\frac{-2\pi}{(\lambda
(z)+\|\hat{z}_{1}\|)^{2}}(\frac{\partial\lambda}{\partial z_{1}}+\frac{z_{1}%
}{\|\hat z_{1}\|}),$ and $\frac{\partial^{2}A(z)}{\partial z_{1}\partial z_{4}}$ is
bounded.\ It is multiplied \ by $f(z,A(z)),$ which tends to zero when $\hat
{z}_{1}$ tends to zero.\ Then it is zero at points $\hat{z}_{1}=0,$ and it is continuous.

Case(2):$\frac{\partial A(z)}{\partial z_{1}}=\frac{-2\pi}{(\lambda
(z)+\|\hat{z}_{1}\|)^{2}}(\frac{\partial\lambda}{\partial z_{1}}+\frac{z_{1}%
}{\|\hat z_{1}\|}),$ and $\frac{\partial^{2}A(z)}{\partial z_{1} \partial z_{2}%
}=C(z)+D(z)\frac{z_{1}z_{2}}{||\hat{z}_{1}||^ {3}},$ where $C(z)$ is
bounded, $D(z)$ is continuous. Then, the question is the continuity to zero of
$\varphi(z)=\frac{f(z,A(z))}{||\hat{z}_{1}||},$ in a neighborhood of the set
$E=\{\hat{z}_{1}=0\}.$ Let us use Formula (\ref{expf}). It gives
$f(z,A(z))=\tilde{Q}_{z,r}(\hat{z}_{1},A(z)-\frac{2\pi}{\lambda(z)}).$ But
$A(z)=\frac{2\pi}{\lambda(z)+||\hat{z}_{1}||},$ then, $A(z)-\frac{2\pi
}{\lambda(z)}=\psi(z)||\hat{z}_{1}||,$ where $\psi(z)$ is continuous. It
follows that $\varphi(z)$ tends to zero when $\hat{z}_{1}$ tends to zero. The
 sub-Riemannian volume is $\mathcal{C}^{2}$ in a neighborhood of $\tilde
{U}_{s_{0}}.$

{It follows that  $f_{\cal S P}(\xi)$ is generically $\mathcal{C}^{2}$ except on a bad set of codimension $8$ in $S^1\times M$, and 
the theorem is proved.  In the case $n=6,$ the bad
set is generically empty in $ S^1\times M$ and property (R) is open dense in
$\mathcal{S}$.}

\appendix

\section{Appendix\label{app}}

\subsection{Pure Quaternions in $so(4)$\label{quat}}

In $so(4)$, it is natural and useful for computations to use quaternionic
notations. Set:%
\[
i=\left(
\begin{array}
[c]{cccc}%
0 & -1 & 0 & 0\\
1 & 0 & 0 & 0\\
0 & 0 & 0 & -1\\
0 & 0 & 1 & 0
\end{array}
\right)  ,\text{ \ \ }j=\left(
\begin{array}
[c]{cccc}%
0 & 0 & -1 & 0\\
0 & 0 & 0 & 1\\
1 & 0 & 0 & 0\\
0 & -1 & 0 & 0
\end{array}
\right)  ,\text{ \ \ }k=\left(
\begin{array}
[c]{cccc}%
0 & 0 & 0 & -1\\
0 & 0 & -1 & 0\\
0 & 1 & 0 & 0\\
1 & 0 & 0 & 0
\end{array}
\right)  ,
\]

\[
\hat{\imath}=\left(
\begin{array}
[c]{cccc}%
0 & -1 & 0 & 0\\
1 & 0 & 0 & 0\\
0 & 0 & 0 & 1\\
0 & 0 & -1 & 0
\end{array}
\right)  ,\text{ \ \ }\hat{\jmath}=\left(
\begin{array}
[c]{cccc}%
0 & 0 & 1 & 0\\
0 & 0 & 0 & 1\\
-1 & 0 & 0 & 0\\
0 & -1 & 0 & 0
\end{array}
\right)  ,\text{ \ \ }\hat{k}=\left(
\begin{array}
[c]{cccc}%
0 & 0 & 0 & -1\\
0 & 0 & 1 & 0\\
0 & -1 & 0 & 0\\
1 & 0 & 0 & 0
\end{array}
\right)  .
\]

The matrices $i,j,k$ (resp.\ $\hat{\imath},\hat{\jmath},\hat{k})$ generate the
so-called pure quaternions (resp. pure skew-quaternions), the space of which
is denoted by $Q$ (resp. $\hat{Q}).$\ The Lie algebra $so(4)=Q\oplus\hat{Q},$
and quaternions commute with skew-quaternions: $[Q,\hat{Q}]=0.$

We endow $so(4)$ with the Hilbert-Schmidt scalar product: $<L_{1},L_{2}%
>=trace(L_{1}^{\prime}L_{2}).$

Then, $i,j,k,\hat{\imath},\hat{\jmath},\hat{k}$ form an orthonormal basis. The
eigenvalues $\omega_{1},\omega_{2}$ of $A=q+\hat{q}$ meet:%
\begin{equation}
-(\omega_{1,2})^{2}=(||q||\pm||\hat{q}||)^{2}. \label{eigquat}%
\end{equation}

As a consequence, an element $A\in so(4)$ has a double eigenvalue iff $A\in
Q\cup\hat{Q}.$

\subsection{Versal deformation of skew-symmetric matrices\label{versal}}

The results of Arnold in \cite{Arnold} can be easily extended to the real
smooth case ($\mathcal{C}^{\infty}),$ for skew-symmetric matrices, under the
action of the orthogonal group:

\begin{theorem}
\cite{Arnold} Let $N(p)$ be a family of $n\times n$ matrices smoothly
depending on $p$ at $(\mathbb{R}^{l},0).$ Let  $O_{N}$ be the orbit of $N=N(0)$
under the action of $Gl(n,\mathbb{R)}$ by conjugation. Let $T(\mu)$ be a
smooth family of matrices, depending on the parameter $\mu\in\mathbb{R}^{k},$
such that the  mapping $\mu\rightarrow T(\mu)$ \textbf{transversally}
intersects $O_{N}$ at some $\tilde{N}=g^{-1}Ng.\ $Then, there is a family of
(smoothly depending on $p)$ matrices $g(p)$ and a smooth mapping
 $p\rightarrow\mu(p),$ such that  $N(p)=g(p)^{-1} T(\mu(p))g(p).$ Moreover, for the transversal $T(\mu),$ one can  choose the
centralizer of $N$ in $gl(n,\mathbb{R)}.$
\end{theorem}

We rephrase the result in the case of a skew-symmetric matrix $N$ that has a
double (but not triple) eigenvalue. Then, by section \ref{quat}, we can assume
that $N$ is (conjugate to) a block-diagonal $Bd(\alpha\hat{q},\delta),$
where $\hat{q}$ is a unit skew-quaternion and $\delta$ is a block-diagonal
skew symmetric matrix with $2\times2$ blocks and non multiple eigenvalues. The
centralizer of $\hat{q}$ in $so(4,\mathbb{R)}$ is the vector space of matrices
of the form $\lambda\hat{q}+q,$ where $q$ varies over pure quaternions. Then,
the centralizer of $N$ in $so(n,\mathbb{R)}$ is the space of block diagonal
matrices Bd$(\lambda\hat{q}+q,\Delta),$ where $q$ varies over pure quaternions
and $\Delta$ varies over $2\times2$ skew-symmetric block diagonal matrices.

Hence, we can find a smooth $g(p)\in SO(n,\mathbb{R)}$, and a smooth
$\mu(p)$ such that:%
\begin{align}
N(p)  &  =g(p)^{-1}T(\mu(p))g(p),\text{ \ \ with}%
\label{versalform}\\
T(\mu)  &  =\mbox{Bd}(\lambda(\mu)\hat{q}+q(\mu),\Delta(\mu)).\nonumber
\end{align}

The versal deformation $T(\mu)$ is not universal (which means that
$\mu(p)$ is not uniquely determined by $N(p)),$ however, the
nondiagonal eigenvalues of $T(\mu)$ are given by the Formula (\ref{eigquat})%
.\ It follows that $q$ is determined modulo conjugation by a unit quaternion.
On the other hand, the functions $\lambda(\mu),\Delta(\mu)$ are smooth and
$\lambda(\mu)$ is nonzero.

\subsection{{Codimension of double and triple eigenvalues}}
\label{a-codimension}
\begin{lemma}
\label{l1}
We have the following: \\[1mm]
(i) the set of skew symmetric matrices with a double eigenvalue is an
algebraic subset of codimension 3 in skew symmetric matrices;\\[.7mm]
(ii) the set of skew symmetric matrices with a triple eigenvalue is an
algebraic subset of codimension 8 in skew symmetric matrices.
\end{lemma}

The proof  of (i) is given in the appendix of \cite{rom}. The proof of (ii) given hereafter is a generalization.
We restrict ourself to the even dimensional case {$so(2n)$}, the odd dimensional case being similar.

We consider the set ${\cal D}$ of block-diagonal matrices $D$ of the form
$$
D=\mbox{Bd}(\alpha J,\alpha J,\alpha J,\alpha_4 J,\ldots,\alpha_n J),
$$
of dimension $N=2n$, with $J=\left(\begin{array}{cc} 0&-1\\1&0 \end{array}     \right)$ and we show only that the union of the orbits under orthogonal conjugation of the elements of ${\cal D}$ has codimension $8$ at least. To do this we consider generic elements of ${\cal D}$ only: for non-generic elements the dimension of the orbit is smaller.

To compute the dimension of the orbit ${\cal O}_D$ of $ D$ it is enough to compute the dimension of the stabilizer $G$ of $D$, and then to compute the dimension of the Lie algebra ${\cal L}=$Lie$(G)$, which is just the centralizer ${\cal C}$ of $D$.

By a direct computation one gets (if $D$ is a generic element) that elements $C$ of ${\cal C}$ are of the form
$$
C=\mbox{Bd}(A_1,\Delta),
$$
where $A_1$ is  $6\times6$ and $\Delta$ is block diagonal with $2\times2$ blocks. Both $A_1$ and $\Delta$ are skew-symmetric and $A_1$ is of the form
$$
\left(\begin{array}{ccc} \alpha_1 J &B_{1,2}&B_{1,3}\\-B_{1,2}& \alpha_2 J& B_{2,3}\\
-B_{1,3}&-{B_{2,3}}&\alpha_3 J 
\end{array}\right),
$$ 
and $$B_{i,j}=\left(\begin{array}{cc} \beta_{i,j}&\gamma_{i,j}\\-\gamma_{i,j}&\beta_{i,j} \end{array}     \right)
$$
Then, dim$(C)=$dim$(G)=n-3+9=n+6$. Therefore dim$({\cal O}_D)=m-n-6$, {where $m=n(2n-1)$ is the dimension of $so(2n)$. }The dimension of ${\cal D}$ is  $n-2$. Hence the dimension of the union of the orbits through points of ${\cal D}$ is $m-n-6+n-2=m-8.$

\subsection{Proof of Lemma \ref{l-genR}: Genericity of (R)}
\label{genericityR}

{We consider the set $\mathcal{S}$ of corank-2  sub-Riemannian metrics on a fixed
manifold $M$, equipped with the Whitney topology. The result being essentially
local, we may assume that $M$ is an open set of $\mathbb{R}^{n},$ with global
coordinates $\xi$,\ and that our  sub-Riemannian metrics are globally specified
by an orthonormal frame, i.e. $s=(F_{1},...,F_{p}).\ $


\medskip
For the moment, we fix $s\in\mathcal{S}.$} We consider two independant one
forms $\omega_{1},\omega_{2}$ on $M,$ that vanish on $\Delta,$ and we set
$\tilde{L}_{i}=d\omega_{i|\Delta},$ and $L_{i}$ is the skew-symmetric matrix
defined by $\tilde{L}_{i}$ via the metric, and moreover   we impose (as in
\ref{popp}) that $L_{1}(\xi),L_{2}(\xi)$ are Hilbert-Schmidt-orthonormal. The
matrices $L_{1},L_{2}$ are defined uniquely modulo a rotation $\hat{L}%
_{1}=\cos(\alpha(\xi))L_{1}+\sin(\alpha(\xi))L_{2},$ $\hat{L}_{2}%
=-\sin(\alpha(\xi))L_{1}+\cos(\alpha(\xi))L_{2}.\ $They are the same as the
matrices $L_{i}$ in Section \ref{nilap} and they meet: $(L_i)_{k,l}=\omega
_{i}([F_{k},F_{l}])=d\omega_{i}(F_{k},F_{l}).$

In coordinates, we set $z=(\xi,\theta),$ and $A(z)=\cos(\theta)L_{1}(\xi
)+\sin(\theta)L_{2}(\xi).$

We fix a point a point $z_{0}=(\theta_{0},\xi_{0})\in\tilde{U}_{s}\subset
 S^1\times M$, and we work in a neighborhood of $z_{0}.\ $By what has just
beeen said, we can perform a constant rotation to have $\theta_{0}=0.\ $

\bigskip\smallskip

Local coordinates $\xi=(x,y)$ in $M$ around $\ \xi_{0}$ can be found, with
$x_{0}=0,$ $y_{0}=0,$ such that:\\[2mm]
1.$\ F_{i}(\xi_{0})=\frac{\partial}{\partial x_{i}},$ $i=1,\ldots,p$ \\[2mm]
2.$\ \omega_{j}(\xi
_{0})=dy_{j}-x^{\prime}{L}_{j}(\xi_{0})dx,$ $j=1,2$ \\[2mm]
 3.$\ A(z_{0})=L_{1}(\xi _{0})$ is $2\times2$ \ block diagonal with decreasingly ordered (moduli of)
eigenvalues. (For this last point, we use an (irrelevant) rotation in the
distribution $\Delta_{\xi_{0}},$ i.e. a constant rotation of the orthonormal frame)

\begin{remark}
Note that at a point $z_{0}=(\xi_{0},\theta_{0})\in\tilde{U}_{s}$, the two
(moduli of) highest eigenvalues of $A(z_{0})$ are equal. However, the whole
construction here holds at each point of $ S^1\times M.$
\end{remark}

\bigskip In these coordinates, we can write (locally) $s$ in the following
form: $F_{i}(\xi)=(e_{i}+B^{i}\xi)\frac{\partial}{\partial x}+x^{\prime}%
L_{1}(\xi_{0})e_{i}\frac{\partial}{\partial y_{1}}+x^{\prime}L_{2}(\xi
_{0})e_{i}\frac{\partial}{\partial y_{2}}+O^{2}(\xi),$ where $e_{i}%
=(0,..,1,..,0)$ is the $i^{th}$ coordinate vector in $\mathbb{R}^{p},$ where
$B^{i}$ is a $p\times n$ matrix, and $O^{2}(\xi)$ is a term of order 2 in
$\xi,$ i.e. $O^{2}(\xi)$ is in $\mathcal{I}^{2},$ where $\mathcal{I}$ is the
ideal of smooth germs at 0 in $\mathbb{R}^{n},$ generated by the components
$\xi_{i}.$

This choice of notations for the vector fields $F_{i}$ is adapted to the
transversality arguments we want to apply later.\ Note that $L_{1}(\xi_{0}),$
$B^{i}$ are, in coordinates, components of the one-jet $j^{1}s(\xi_{0})$ of
$s$ at $\xi_{0}.$

\bigskip

Define the $p\times p$ matrix $U^{r}$ by $U_{i,j}^{r}=B_{i,r}^{j},$ and by
$\Omega_{1},\Omega_{2}$ the skew-symmetric matrices associated with the
2-forms $d\omega_{1|\Delta}(\xi_{0}),d\omega_{2|\Delta}(\xi_{0})$ in the
chosen coordinates. It is not hard to compute the tangent mappings $TL_{1}%
(\xi_{0})$ and $TL_{2}(\xi_{0})$ (we temporarily write $TL(\xi_{0})$ and
$\Omega$ for convenience):
\begin{equation}
TL(\xi_{0})(e_{r})=U^{r\prime}\Omega-\Omega^{\prime}U^{r},\text{
\ \ }r=1,...,n.\label{calT}%
\end{equation}

For this, one just uses $d\circ d=0,$ and $TL_{k,l}(\xi_{0})(e_{r}%
)=d\omega(TF_{k}(e_{r}),F_{l})+d\omega(F_{k},TF_{l}(e_{r})),$ where $d\omega$
stands for $d\omega_{1}$ or $d\omega_{2}.$

\bigskip\bigskip On the other hand, we have, using the versality theorem in a
neighborhood of $z_{0}$:%
\begin{equation}
A(z)=\cos(\theta)L_{1}(\xi)+\sin(\theta)L_{2}(\xi)=H(z)\mbox{Bd}(\lambda(z)\hat
{q}+q(z),\Delta(z))H^{\prime}(z), \label{dec}%
\end{equation}

 in which we already assumed that $\theta_{0}=0,$ and the coordinates
$\xi=(x,y)$ were already chosen for $A(z_{0})$ to be diagonal. Also,
$H(z_{0})=Id.$

\begin{remark}
1. The decomposition (\ref{dec}) is not unique: the quaternion $q(z)$ is
defined modulo conjugation by a unit quaternion, $\tilde{q}(z)=q_{1}%
(z)q(z)q_{1}(z)^{-1}.$ However, the tangent mapping $Tq(z)$ is changed for
$T\tilde{q}(z)=[Tq_{1}(z),q(z)]+q_{1}(z)Tq(z)q_{1}(z)^{-1}.$ But on $\tilde
{U}_{s},$ $q(0)=0,$ hence the rank of $Tq(z_{0})$ remains unchanged.

2.\ The decomposition can easily be made unique, by making (following Arnold
\cite{Arnold}) some particular choice of a (mini)transversal to the
centraliser of $A(z_{0})$.\ For instance, one could chose the
(Hilbert-Schmidt) orthogonal supplement to the centralizer {of $A(z_{0})$}
through $A(z_{0}).$
\end{remark}

Let $\Pi_{Q}$ $:so(n)\rightarrow Q\simeq\mathbb{R}^{3},$ be the projection
associating to the matrices, the quaternionic components of the first
$4\times4$ block on the diagonal.

By (\ref{dec}), we have:%
\begin{equation}
Tq(z_{0})(V)=\Pi_{{Q}}TA(z_{0})(V)+\Pi_{{Q}}[TH^{\prime}(z_{0})(V),\mbox{Bd}(z_{0})].
\label{Tq}%
\end{equation}
We can consider the fiber mapping $\pi_{Q}:J^{1}\mathcal{S\times}%
S^1\rightarrow Q\times\mathcal{M}(3,n+1),$ $\pi_{Q,z_{0}}:(j^{1}%
s\mathcal{(}\xi_{0}),\theta_{0})\rightarrow(q(z_{0}),Tq(z_{0}))$ 
($\mathcal{M}(3,n+1)$ being the set of 
$3\times(n+1)$ real matrices),%
\[
\pi_{Q,z_{0}}(L_{1},L_{2},B^{i},i=1,...,r)=\{\Pi_{Q}(A(z_{0})),\text{ }\Pi
_{Q}\circ TA(z_{0})\}.
\]
The following lemma is an easy consequence of (\ref{calT}), (even easier to
prove if one considers that\footnote{Note that in fact, in the chosen
coordinates, $\Omega_{i}=L_{i}$ since F$_{i}(z_{0})=e_{i}$ on $\Delta_{z_{0}}%
$} $\Omega=\Omega_{1}=L_{1},$ is $2\times2$ block diagonal, the 2 first blocks
being both nonzero):

\begin{lemma}
\label{sub}The mapping $\pi_{Q,z_{0}}$ is a linear submersion.
\end{lemma}

It follows from Lemma \ref{sub} that the mapping $\rho:J^{1}\mathcal{S\times
}S^1\rightarrow\mathbb{R}^{3}\times\mathcal{M}(3,n+1),$ $(z_{0},L_{1}%
,L_{2},B^{i},i=1,...,{p})\rightarrow Tq(z_{0})$ is a submersion.

The codimension $d_{0\text{ }}$of the algebraic set of $3\times(n+1)$ matrices
that have corank 1 at least is $d_{0}=(n-1)$ [product of coranks in the
$3\times(n+1)$ matrices]. By Lemma \ref{l1}, the set of skew-symmetric
matrices that have double maximum eigenvalue is $d_{1}=3$. Therefore, by the
transversality theorems \cite{abra}, there is a residual subset  of the set of 
 sub-Riemannian metrics, for which the codimension of the set of $z=(\theta,\xi)$
in $S^1\times M$ where $A(z)$ corresponds to a double eigenvalue, an
property (R) holds at $(\theta,\xi)$, is a stratified set of codimension
$d_{0}+d_{1}=n+2.$

\subsection{Volume of the unit ball}
\label{formula}
We keep the notations of Section \ref{proof}.

\[
\mathcal{E}_{t}\mathcal{(}p^x_{0},p^y_{0})= \mathcal{E}_{t}\mathcal{(}p^x_{0},r,\theta)= \mathcal{E}_{1}\mathcal{(}t \, p^x_{0},t \,r,\theta),
\]
ant $ t_{cut}=A_\xi(\theta)/r$. The domain of $\mathcal{E}_{t}$ for the unit ball is 
$$
D_{\mathcal{E}_{t}}=\left\{ (p^x_0,r,\theta,t)~ |
~ \theta\in[0,2\pi], ~p^x_0\in S^{{p-1}}, t\in\min(1,t_{cut}),~r\in[0,+\infty[
\right\},
$$
where $S^{p-1}$ denotes the unit Euclidean sphere in $\mathbb{R}^{p}$ and later $B$ denotes the unit Euclidean ball $\mathbb{R}^{p}$. 
In this formula, the boundary of this set in the variables $r,t$ is parametrized by $r$. Equivalently if we parametrize this boundary by $t$ we get,
$$
D_{\mathcal{E}_{t}}=\left\{ (p^x_0,r,\theta,t)~ |
~ \theta\in[0,2\pi], ~p^x_0\in S^{p-1}, t\in[0,1 ],~r=A_\xi(\theta)/t
\right\}.
$$
Now set $\tilde p=t\,p^x_0,$ $\tilde r=t\,r$. This implies  for the domain $D_{\mathcal{E}_{1}}$ of $\mathcal{E}_{1}(\tilde p ,\tilde r, \theta),$

$$D_{\mathcal{E}_{1}}=\left\{ (\tilde p, \tilde r,\theta)~ |
~ \theta\in[0,2\pi], ~\tilde p\in B,~\tilde r\in[0,A_\xi(\theta)]
\right\}.
$$
The volume of the unit ball of the nilpotent approximation at $\xi$ is 

$$
V_\xi=\int_{\mathcal{E}_{1}(D_{\mathcal{E}_{1}})}\mbox{Popp}=\int_{\mathcal{E}_{1}(D_{\mathcal{E}_{1}})}dx\wedge dy,
$$
that is 
$$
V_\xi=\int_0^{2\pi} \int_0^{A_\xi(\theta)}\int_{B} J_{\mathcal{E}_{1}(\tilde p.\tilde r,\theta)d\tilde p\,d\tilde r\,d\theta},
$$
where  $J_{\mathcal{E}_{1}}$ is the Jacobian determinant of $\mathcal{E}_{1}(\tilde p ,\tilde r, \theta)$.

\begin{acknowledgement}
we thank with great respect the memory of\ V.\ Zakalyukin who gave the basic
idea of this result, in june 2011.

 This research has been supported by the European Research Council, ERC StG 2009 ``GeCoMethods'', contract number 239748, by the ANR Project GCM, program ``Blanche'', project number NT09-504490 and by the DIGITEO project CONGEO.

\end{acknowledgement}

\end{document}